\documentclass[12pt]{article}

\renewcommand{\leq}{\leqslant}
\renewcommand{\geq}{\geqslant}

\usepackage{amsthm,amsmath,amssymb,latexsym}

\usepackage{epsfig,cite,lineno}

\newtheorem{thm}{Theorem}[section]
\newtheorem{prop}[thm]{Proposition}

\newtheorem{cor}[thm]{Corollary}

\newtheorem{conj}[thm]{Conjecture}

\theoremstyle{remark}

\numberwithin{equation}{section}

\begin{document}

\begin{center}
{\Large\bf Factors of sums involving $q$-binomial coefficients\\[5pt] and powers of $q$-integers}
\end{center}
\vskip 2mm \centerline{Victor J. W. Guo$^1$ and Su-Dan Wang$^{2}$}
\begin{center}
{\footnotesize $^1$School of Mathematical Sciences, Huaiyin Normal University, Huai'an, Jiangsu 223300,\\
 People's Republic of China\\
{\tt jwguo@hytc.edu.cn}\\[10pt]
%$^2$Universit\'e de Lyon, Lyon, F-69003, France\\

$^2$Department of Mathematics, East China Normal University, Shanghai 200062,\\
 People's Republic of China\\
 {\tt sudan199219@126.com}
 }
\end{center}

%\centerline{September 3, 2009}

\vskip 0.7cm {\noindent{\bf Abstract.}
We show that, for all positive integers
$n_1, \ldots, n_m$, $n_{m+1}=n_1$, and any non-negative integers $j$ and $r$ with $j\leqslant m$, the expression
$$
\frac{1}{[n_1]}{n_1+n_{m}\brack n_1}^{-1}
\sum_{k=1}^{n_1}[2k][k]^{2r}q^{jk^2-(r+1)k}\prod_{i=1}^{m} {n_i+n_{i+1}\brack n_i+k}
$$
is a Laurent polynomial in $q$ with integer coefficients, where $[n]=1+q+\cdots+q^{n-1}$ and ${n\brack k}=\prod_{i=1}^k(1-q^{n-i+1})/(1-q^i)$.
This gives a $q$-analogue of a divisibility result on the Catalan triangle obtained by the first author and Zeng, and
also confirms a conjecture of the first author and Zeng. We further propose several related conjectures.
}



%%%%%%%%%%%%%%%%%%%%%%%%%%%%%%%%%%%%%%%%%%%%%%%%%%%%%%%%%%%%%%%%%%%%%%%%%%%%%%%%%%%%%%%%%%%%%%%%%%%
\section{Introduction}
A few years ago, motivated by the divisibility results in \cite{Calkin,CD,CC,GJZ,GHMR,MR,Shapiro,Zu},  the first author and Zeng \cite{GZ2010} proved that
\begin{align}
\frac{2}{n^2}{2n\choose n}^{-1}\sum_{k=1}^n{2n\choose n-k}k^{2r+1}\in\mathbb{Z}\quad\text{for } r\geqslant 1,\label{eq:guo1}
\end{align}
and more general
\begin{align}
\frac{2}{n_1}{n_1+n_{m}\choose n_1}^{-1}
\sum_{k=1}^{n_1}k^{2r+1}\prod_{i=1}^{m} {n_i+n_{i+1}\choose n_i+k}\in\mathbb{Z} \quad\text{for } r\geqslant 0, \label{eq:guo2}
\end{align}
where $n_{m+1}=n_1$.

The first objective of this paper is to give a $q$-analogue of \eqref{eq:guo1} and \eqref{eq:guo2}. Recall that
the {\it $q$-integers} are defined by $[n]=\frac{1-q^n}{1-q}$ and the {\it $q$-shifted factorials} (see \cite{GR}) are defined by $(a;q)_0=1$
and $(a;q)_n=(1-a)(1-aq)\cdots (1-aq^{n-1})$ for $n=1,2,\ldots.$ The {\it $q$-binomial coefficients} are defined as
$$
{n\brack k}=
\begin{cases}\displaystyle\frac{(q;q)_n}{(q;q)_k (q;q)_{n-k}}, &\text{if $0\leqslant k\leqslant n$,} \\[10pt]
0, &\text{otherwise.}
\end{cases}
$$
Our main results can be stated as follows.
\begin{thm}\label{thm:half}
For $j=0,1$ and all positive integers $n$ and $r$, the expression
\begin{align}
\frac{1}{[n]^2}{2n\brack n}^{-1}
\sum_{k=1}^{n}[2k][k]^{2r}q^{(r+1)(n-k)+jk^2}{2n\brack n+k} \label{eq:main-one}
\end{align}
is a polynomial in $q$ with integer coefficients.
\end{thm}

\begin{thm}\label{thm:nnnhalf}
Let $n_1,\ldots,n_{m},n_{m+1}=n_1$ be positive integers. Then for any non-negative integers $j$ and $r$ with $j\leqslant m$, the expression
\begin{align}
\frac{1}{[n_1]}{n_1+n_{m}\brack n_1}^{-1}
\sum_{k=1}^{n_1}[2k][k]^{2r}q^{jk^2-(r+1)k}\prod_{i=1}^{m} {n_i+n_{i+1}\brack n_i+k}   \label{eq:main-two}
\end{align}
is a Laurent polynomial in $q$ with integer coefficients.
\end{thm}

One may wonder why to consider $[2k][k]^{2r}/[2]$ as a $q$-analogue of $k^{2r+1}$ rather than $[k]^{2r+1}$.
The idea comes from Schlosser's work \cite{Schlosser} on $q$-analogues of sums of powers of the first $n$ consecutive natural numbers and
the first author and Zeng's work \cite{GZ04} which positively answers Schlosser's question.

The first author and Zeng \cite{GZ2010} also  introduced the  \emph{$q$-Catalan triangle} with entries given by
$$
B_{n,k}(q):=\frac{[k]}{[n]}{2n\brack n-k},\quad1\leqslant k\leqslant n.
$$
The second objective is to give the following congruence related to the $q$-Catalan triangle.
\begin{thm}\label{thm:ppower}
Let $n$ be a positive integer. Let $r\geq 0$ and $s\geq 1$ such that $r\not\equiv s\pmod 2$. Then, for $0\leqslant j\leqslant s$, the expression
\begin{align*}
{2n\brack n}^{-1}\sum_{k=1}^{n}(1+q^k)[k]^{r}B_{n,k}^{s}(q)q^{jk^2-(r+s+1)k/2}
\end{align*}
is a Laurent polynomial in $q$ with integer coefficients.
\end{thm}

It is easy to see that Theorem \ref{thm:ppower} is a generalization of \cite[Theorem 1.4]{GZ2010}.
Letting $q=1$ in Theorem \ref{thm:ppower}, we confirm a conjecture in \cite{GZ2010}.

The paper is organized as follows. We shall prove Theorems \ref{thm:half}--\ref{thm:ppower}
in Sections 2--4, respectively. The $q$-Pfaff-Saalsch\"utz identity will play an important role in our proof.
We give some consequences of Theorem \ref{thm:half} in Section 5. Finally, some open
problems will be proposed in Section 6.

%%%%%%%%%%%%%%%%%%%%%%%%%%%%%%%%%%%%%%%%%%%%%%%%%%%%%%%%%%%%%%%%%%%%%%%%%%%%%%%%%%%%%%%%%%%%%%%%%%%%%%%%%%%%%%%%%%%%%
\section{Proof of Theorem \ref{thm:half}}
We first consider the $j=0$ case. Let
\begin{align*}
S_r(n;q)=\sum_{k=1}^{n}[2k][k]^{2r}q^{(r+1)(n-k)}{2n\brack n+k}.
\end{align*}
It is easy to see that
\begin{align*}
S_0(n;q)=\sum_{k=1}^n[2k]q^{n-k}{2n\brack n-k}
&=\sum_{k=1}^n [n-k+1]{2n\brack n-k+1}
-\sum_{k=1}^n [n-k]{2n\brack n-k}\\
&=[n]{2n\brack n}.
\end{align*}

For $r\geq1$, noticing the relation
\begin{align*}
[k]^2{2n\brack n+k}q^{n-k}=[n]^2{2n\brack n+k}-[2n][2n-1]{2n-2\brack n+k-1},
\end{align*}
we have
\begin{align}
S_r(n;q)=[n]^2S_{r-1}(n;q)-q^r[2n][2n-1] S_{r-1}(n-1;q),\ n=1,2,\ldots. \label{eq:sumnkl}
\end{align}
It follows that
\begin{align*}
S_1(n;q)=[n]^3{2n\brack n}-q[2n][2n-1][n-1]{2n-2\brack n-1}=[n]^2{2n\brack n}.
\end{align*}
This proves that $S_1(n;q)\equiv 0\pmod{[n]^2{2n\brack n}}$. Applying the recurrence relation \eqref{eq:sumnkl} and by induction on $r$,
we can prove that, for all positive integers $r$, there holds
\begin{align}
S_r(n;q)\equiv 0\pmod{[n]^2{2n\brack n}}. \label{eq:srnq}
\end{align}

For $j=1$, let
\begin{align*}
T_r(n;q)=\sum_{k=1}^{n}[2k][k]^{2r}q^{(r+1)(n-k)+k^2}{2n\brack n+k}.
\end{align*}
It is well known that ${n\brack k}_{q^{-1}}={n\brack k}q^{k(k-n)}$. Hence,
$ T_r(n;q)=q^{n^2+2rn+2n-2r-1} S_r(n;q^{-1})$. By \eqref{eq:srnq}, we have
\begin{align*}
T_r(n;q)\equiv 0\pmod{[n]^2{2n\brack n}}.
\end{align*}
Since both $S_r(n;q)$ and $T_r(n;q)$ are polynomials in $q$, we obtain the desired conclusion.

%%%%%%%%%%%%%%%%%%%%%%%%%%%%%%%%%%%%%%%%%%%%%%%%%%%%%%%%%%%%%%%%%%%%%%%%%%%%%%%%%%%%%%%%%%%%%%%%%%%%%%%%%%%%%%%%%%%%%
\section{Proof of Theorem \ref{thm:nnnhalf}}
We will need  the $q$-Pfaff-Saalsch\"utz identity (see \cite[(II.12)]{GR} or \cite{GJZ2005,Zeilberger}):
\begin{align}
{n_1+n_2\brack n_1+k}{n_2+n_3\brack n_2+k}{n_3+n_1\brack n_3+k}
=\sum_{s=0}^{n_1-k}\frac{q^{s^2+2ks} (q;q)_{n_1+n_2+n_3-k-s}}
{(q;q)_s (q;q)_{s+2k}(q;q)_{n_1-k-s}(q;q)_{n_2-k-s}(q;q)_{n_3-k-s}}, \label{eq:qpfaff}
\end{align}
where $\frac{1}{(q;q)_n}=0$ if $n<0$. Denote \eqref{eq:main-two} by $S(n_1,\ldots,n_{m};r,j,q)$.
Namely,
\begin{equation}\label{eq:rewriting}
S(n_1,\ldots,n_{m};r,j,q)
=\frac{(q;q)_{n_1-1} (q;q)_{n_{m}}}{(q;q)_{n_1+n_{m}}}
\sum_{k=1}^{n_1} [2k][k]^{2r}q^{jk^2-(r+1)k} C(n_1,\ldots,n_{m};k),
\end{equation}
where
$$
C(a_1,\ldots,a_l;k)=\prod_{i=1}^l {a_i+a_{i+1}\brack a_i+k}\quad (a_{l+1}=a_1).
$$

It is easy to see that, for $m\geq 3$,
\begin{align}
C(n_1,\ldots,n_m;k)=\frac{(q;q)_{n_2+n_3}(q;q)_{n_m+n_1}}{(q;q)_{n_1+n_2}(q;q)_{n_m+n_3}}
{n_1+n_2\brack n_1+k}{n_1+n_2\brack n_2+k}C(n_3,\ldots,n_m;k).  \label{eq:cn1nk}
\end{align}
Letting $n_3\to\infty$ in \eqref{eq:qpfaff}, we get
\begin{align}
{n_1+n_2\brack n_1+k}{n_1+n_2\brack n_2+k}
=\sum_{s=0}^{n_1-k}\frac{q^{s^2+2ks}(q;q)_{n_1+n_2}}{(q;q)_{s}(q;q)_{s+2k}(q;q)_{n_1-k-s}(q;q)_{n_2-k-s}}.\label{eq:qpfaffn1n2}
\end{align}
Substituting \eqref{eq:cn1nk} and \eqref{eq:qpfaffn1n2} into \eqref{eq:rewriting}, we obtain
\begin{align*}
&\hskip -2mm S(n_1,\ldots,n_{m};r,j,q) \\
&=\frac{(q;q)_{n_2+n_3}(q;q)_{n_1-1} (q;q)_{n_m}}{(q;q)_{n_m+n_3}}\sum_{k=1}^{n_1}\sum_{s=0}^{n_1-k}
\frac{q^{s^2+2ks+jk^2-(r+1)k}[2k][k]^{2r} C(n_3,\ldots,n_{m};k)}{(q;q)_{s}(q;q)_{s+2k}(q;q)_{n_1-k-s}(q;q)_{n_2-k-s}}  \\
&=\frac{(q;q)_{n_2+n_3}(q;q)_{n_1-1} (q;q)_{n_m}}{(q;q)_{n_m+n_3}}\sum_{l=1}^{n_1} q^{l^2}\sum_{k=1}^{l}
\frac{q^{(j-1)k^2-(r+1)k}[2k][k]^{2r} C(n_3,\ldots,n_{m};k)}{(q;q)_{l-k}(q;q)_{l+k}(q;q)_{n_1-l}(q;q)_{n_2-l}},
\end{align*}
where $l=s+k$. Observing that, for $m\geqslant 3$,
$$
\frac{C(n_3,\ldots, n_{m};k)}{(q;q)_{l-k}(q;q)_{l+k}}
=\frac{(q;q)_{n_{m}+n_3}}{(q;q)_{n_3+l}(q;q)_{n_{m}+l}}C(l,n_3,\ldots, n_{m};k),
$$
we get the following recurrence relation
\begin{align}
S(n_1,\ldots,n_{m};r,j,q)
=\sum_{l=1}^{n_1} q^{l^2}{n_1-1\brack l-1}{n_2+n_3\brack n_2-l} S(l,n_3,\ldots,n_{m};r,j-1,q). \label{eq:recsr}
\end{align}
On the other hand, for $m=2$, applying \eqref{eq:qpfaffn1n2} we may deduce that
\begin{align}\label{eq:sn1n2}
S(n_1,n_{2};r,j,q)=\sum_{l=1}^{n_1} q^{l^2}{n_1-1\brack l-1}{n_2\brack l}S(l;r,j-1,q).
\end{align}

We now proceed by induction on $m$.  For $m=1$, the conclusion follows readily from the proof of Theorem \ref{thm:half}.
Suppose that the expression $S(n_1,\ldots,n_{m-1};r,j,q)$ is a Laurent polynomial in $q$ with integer coefficients for some $m\geqslant 2$ and $0\leq j\leq m-1$.
Then by the recurrence \eqref{eq:recsr} or \eqref{eq:sn1n2}, so is $S(n_1,\ldots,n_{m};r,j,q)$ for  $1\leq j\leq m$. It is easy to see that
\begin{align*}
S(n_1,\ldots,n_{m};r,0,q)=S(n_1,\ldots,n_{m};r,m,q^{-1}) q^{n_1n_2+n_2n_3+\cdots+n_{m-1}n_m-n_1-2r}\quad\text{for }m\geqslant 2.
\end{align*}
Therefore, the expression $S(n_1,\ldots,n_{m};r,0,q)$ is also a Laurent polynomial in $q$ with integer coefficients.
This completes the inductive step of the proof.

%%%%%%%%%%%%%%%%%%%%%%%%%%%%%%%%%%%%%%%%%%%%%%%%%%%%%%%%%%%%%%%%%%%%%%%%%%%%%%%%%%%%%%%%%%%%%%%%%%%%%%%%%%%%%%%%%%%%%
\section{Proof of Theorem \ref{thm:ppower}}

Let $\Phi_n(x)$ be the $n$-th {\it cyclotomic polynomial}. We need the following result (see \cite[Equation (10)]{KW} or \cite{CH,GZ06}).
\begin{prop}\label{prop:factor}
The $q$-binomial coefficient ${m\brack k}$ can be factorized into
$$
{m\brack k}=\prod_{d}\Phi_d(q),
$$
where the product is over all positive integers $d\leq m$ such that
$\lfloor k/d\rfloor+\lfloor (m-k)/d\rfloor<\lfloor m/d\rfloor$.
\end{prop}

By Proposition \ref{prop:factor}, for any positive integer $n$, we have
\begin{align}
\gcd\left({2n\brack n},[n]\right)=1.  \label{eq:2nbrackn}
\end{align}
Let $r+s\equiv 1\pmod 2$, $m=s\geq 1$, and $0\leqslant j\leqslant s$. Setting $n_1=\cdots=n_s=n$ in Theorem~\ref{thm:nnnhalf}, we see that
\begin{align*}
\frac{1}{[n]}{2n\brack n}^{-1}
\sum_{k=1}^{n}[2k][k]^{r+s-1}q^{jk^2-(r+s+1)k/2} {2n\brack n+k}^s
\end{align*}
is a Laurent polynomial in $q$ with integer coefficients.  Note that $B_{n,k}(q)$ is a polynomial in $q$ with integer coefficients (see \cite{GZ2010,GK}).
Hence, $[k]{2n\brack n+k}=[n]B_{n,k}(q)$ is clearly divisible by $[n]$. It follows that
\begin{align*}
\frac{\gcd\left({2n\brack n},[n]^s\right)}{{2n\brack n}}\sum_{k=1}^{n}(1+q^k)[k]^{r}B_{n,k}^{s}(q)q^{jk^2-(r+s+1)k/2}
\end{align*}
is a Laurent polynomial in $q$ with integer coefficients. The proof then follows from  \eqref{eq:2nbrackn}.

%%%%%%%%%%%%%%%%%%%%%%%%%%%%%%%%%%%%%%%%%%%%%%%%%%%%%%%%%%%%%%%%%%%%%%%%%%%%%%%%%%%%%%%%%%%%%%%%%%%%%%%%%%%%%%%%%%%%%%%%
\section{Some consequences of Theorem \ref{thm:nnnhalf}}
In this section, we will give some consequences of Theorem \ref{thm:nnnhalf} and confirm some conjectures in \cite[Section 7]{GZ2010}.
Letting $n_{2i-1}=m$ and $n_{2i}=n$ for $i=1,\ldots,r$ in Theorem~\ref{thm:nnnhalf} and observing the symmetry of $m$ and $n$,
we obtain
\begin{cor}\label{cor:mnrs}
Let $m$, $n$, and $r$ be positive integers, and let $a$ and $j$ be non-negative integers with $j\leqslant 2r$.
 Then the expression
\begin{align*}
\frac{\gcd([m],[n])}{[m][n]}{m+n\brack m}^{-1}
\sum_{k=1}^{m}[2k][k]^{2a}q^{jk^2-(a+1)k} {m+n\brack m+k}^r {m+n\brack n+k}^r
\end{align*}
is a Laurent polynomial in $q$ with integer coefficients.
\end{cor}

Letting $n_{3i-2}=l$, $n_{3i-1}=m$ and $n_{3i}=n$ for $i=1,\ldots,r$ in Theorem~\ref{thm:nnnhalf}, we get
\begin{cor}
Let $l$, $m$, $n$ and $r$ be positive integers, and let $a$ and $j$ be non-negative integers with $j\leqslant 3r$. Then the expression
\begin{align*}
\frac{\gcd([m],[n])}{[m][n]}{m+n\brack m}^{-1}
\sum_{k=1}^{m}[2k][k]^{2a}q^{jk^2-(a+1)k} {l+m\brack l+k}^r  {m+n\brack m+k}^r {n+l\brack n+k}^r
\end{align*}
is a Laurent polynomial in $q$ with integer coefficients.
\end{cor}

Taking $m=2r+s$ and letting $n_i=n+1$ if $i=1,3,\ldots,2r-1$ and $n_i=n$ otherwise in
Theorem~\ref{thm:nnnhalf}, we get
\begin{cor}
Let $n$, $r$ and $s$ be positive integers, and let $a$ and $j$ be non-negative integers with $j\leqslant 2r+s$. Then the expression
\begin{align*}
\frac{1}{[n][n+1]}{2n+1\brack n}^{-1}
\sum_{k=1}^{n}[2k][k]^{2a}q^{jk^2-(a+1)k} {2n+1\brack n+k+1}^r{2n+1\brack n+k}^r{2n\brack n+k}^s
\end{align*}
is a Laurent polynomial in $q$ with integer coefficients.
\end{cor}

Let $[n]!=[n][n-1]\cdots[1]$. It is clear that Theorem~\ref{thm:nnnhalf} can be restated as follows.

\begin{thm}\label{thm:rennnhalf}
Let $n_1,\ldots,n_{m}$ be positive integers. Let $j$ and $r$ be non-negative integers  with $j\leqslant m$.
Then the expression
\begin{align*}
[n_1-1]!\prod_{i=1}^m\frac{[n_i+n_{i+1}]!}{[2n_i]!}
\sum_{k=1}^{n_1}[2k][k]^{2r}q^{jk^2-(r+1)k}\prod_{i=1}^{m} {2n_i\brack n_i+k},
\end{align*}
where $n_{m+1}=0$, is a Laurent polynomial in $q$ with integer coefficients.
\end{thm}

Letting $n_1=\cdots= n_r=m$ and
$n_{r+1}=\cdots=n_{r+s}=n$ in Theorem~\ref{thm:rennnhalf} and noticing the symmetry of $m$ and $n$, we have

\begin{cor}\label{cor:mn}
Let $m$, $n$, $r$ and $s$ be positive integers, and let $a$ and $j$ be non-negative integers with $j\leqslant r+s$. Then the expression
\begin{align*}
\frac{\gcd([m],[n])[m+n]![m-1]![n-1]!}{[2m]![2n]!}
\sum_{k=1}^{m}[2k][k]^{2a}q^{jk^2-(a+1)k} {2m\brack m+k}^r {2n\brack n+k}^s
\end{align*}
is a Laurent polynomial in $q$ with integer coefficients.
\end{cor}

In particular, the expression
\begin{align*}
\frac{1}{[2n]}{4n\brack n}^{-1}
\sum_{k=1}^{n}[2k][k]^{2a}q^{jk^2-(a+1)k} {4n\brack 2n+k}^r {2n\brack n+k}^s
\end{align*}
is a Laurent polynomial in $q$ with integer coefficients.

The {\it $q$-Narayana numbers} $N_q(n,k)$ and the {\it $q$-Catalan numbers} $C_n(q)$ may be defined as follows:
$$
N_q(n,k)=\frac{1}{[n]}{n\brack k}{n\brack k-1},\quad C_n(q)=\frac{1}{[n+1]}{2n\brack n}.
$$
It is well known that both $q$-Narayana numbers and $q$-Catalan numbers are polynomials in $q$
with non-negative integer coefficients (see \cite{Branden,FH}). It should be mentioned that the definition of $N_q(n,k)$ here differs by a
factor $q^{k(k-1)}$ from that in \cite{Branden}. Motivated mainly by Z.-W. Sun's work on congruences for
combinatorial numbers \cite{Sun0,Sun1,Sun2,Sun3,Sun4,Sun5}, the first author and Jiang \cite{GJ} proved that, for $0\leqslant j\leqslant 2r-1$,
\begin{align}
\sum_{k=-n}^{n}(-1)^{k}q^{jk^2+{k\choose 2}}N_q(2n+1,n+k+1)^r \equiv 0 \pmod{C_n(q)}.
\label{eq:main-2}
\end{align}

Exactly similarly to the proof of \eqref{eq:main-2} in \cite{GJ}, we can deduce the following result from Theorem~\ref{thm:rennnhalf}
by considering four special cases
$(n_1,\ldots,n_{2r})=(n,\ldots,n,n+1,\ldots,n+1),(n+1,\ldots,n+1,n,\ldots,n),(n,n+1,n,n+1,\ldots,n,n+1),(n+1,n,n+1,n,\ldots,n+1,n)$
and noticing that $[n]$ and $[n+1]$ are relatively prime.
\begin{cor}\label{cor:narayana}
Let $n$ and $r$ be positive integers, and let $a$ and $j$ be non-negative integers with $j\leqslant 2r$. Then the expression
\begin{align}
\frac{1}{[n]}{2n\brack n}^{-1}
\sum_{k=1}^{n}[2k][k]^{2a}q^{jk^2-(a+1)k} N_q(2n+1,n+k+1)^r  \label{eq:narayana}
\end{align}
is a Laurent polynomial in $q$ with integer coefficients.
\end{cor}

%%%%%%%%%%%%%%%%%%%%%%%%%%%%%%%%%%%%%%%%%%%%%%%%%%%%%%%%%%%%%%%%%%%%%%%%%%%%%%%%%%%%%%%%%%%%%%%%%%%%%%%%%%%%%%%%%%%%%%%%%
\section{Some open problems}
In this section we propose several related conjectures for further study. Note that some similar conjectures were raised in \cite[Section 3]{GJ}.

It is easy to see that the {\it $q$-super Catalan numbers} $\frac{[2m]![2n]!}{[m+n]![m]![n]!}$ are polynomials in $q$ with integer coefficients.
Warnaar and Zudilin \cite[Proposition 2]{WZ} proved that the $q$-super Catalan numbers are in fact polynomials in $q$
with non-negative integer coefficients. The following conjecture related to the $q$-super Catalan numbers is a refinement of Corollary~\ref{cor:mnrs}
and is also a $q$-analogue of \cite[Conjecture 7.5]{GZ2010}.

\begin{conj}
Let $m$, $n$, $s$ and $t$ be positive integers. Let $r$ be a non-negative integer with $r+s+t\equiv 1\pmod{2}$ and let $j$ be a integer.
Then the expression
\begin{align*}
\frac{[m+n]![m]![n]!}{[2m]![2n]!}
\sum_{k=1}^{m}(1+q^k)[k]^{r}q^{jk^2-(r+s+t+1)k/2} B_{m,k}^s(q) B_{n,k}^t(q)
\end{align*}
is a Laurent polynomial in $q$, and  is a Laurent polynomial in $q$ with non-negative
integer coefficients for $0\leqslant j\leqslant s+t$.
\end{conj}

It seems that Corollary {\rm\ref{cor:narayana}} can be further generalized as follows.
\begin{conj}
Corollary {\rm\ref{cor:narayana}} is still true for any integer $j$, and \eqref{eq:narayana} is a Laurent polynomial in $q$ with non-negative
integer coefficients for $0\leqslant j\leqslant 2r$.
\end{conj}

The following is a generalization of Theorem~\ref{thm:nnnhalf}.

\begin{conj}
Let $n_1,\ldots,n_{m},n_{m+1}=n_1$ be positive integers. Then for any integer $j$ and  non-negative integer $r$, the expression
\begin{align*}
\frac{1}{[n_1]}{n_1+n_{m}\brack n_1}^{-1}
\sum_{k=1}^{n_1}[2k][k]^{2r}q^{jk^2-(r+1)k}\prod_{i=1}^{m} {n_i+n_{i+1}\brack n_i+k}
\end{align*}
is a Laurent polynomial in $q$, and is a Laurent polynomial in $q$ with non-negative integer coefficients
if $0\leqslant j\leqslant m$.
\end{conj}

We end the paper with the following conjecture. Note that when all the $n_i$'s are equal to $n$, it reduces to Corollary {\rm\ref{cor:narayana}}.
\begin{conj}\label{conj:gen-narayana}
Let $n_1,\ldots,n_{m},n_{m+1}=n_1$ be positive integers. Then for any integer $j$ and  non-negative integer $r$, the expression
$$
\frac{1}{[n_1]{n_1+n_m\brack n_1}}\prod_{i=1}^{m}\frac{1}{[n_i+n_{i+1}+1]}\sum_{k=1}^{n_1} [2k][k]^{2r}q^{jk^2-(r+1)k}
\prod_{i=1}^m {n_i+n_{i+1}+1\brack n_i+k}{n_i+n_{i+1}+1\brack n_i+k+1}
$$
is a Laurent polynomial in $q$, and is a Laurent polynomial in $q$ with non-negative integer coefficients
if $0\leqslant j\leqslant 2m$.
\end{conj}

\medskip
\vskip 5mm \noindent{\bf Acknowledgments.} The first author was partially
supported by the National Natural Science Foundation of China (grant 11371144),
the Natural Science Foundation of Jiangsu Province (grant BK20161304),
and the Qing Lan Project of Education Committee of Jiangsu Province.

\renewcommand{\baselinestretch}{1}


\begin{thebibliography}{99}
\small \setlength{\itemsep}{-.8mm}

\bibitem{Branden}P. Br\"and\'en, $q$-Narayana numbers and the flag $h$-vector of $J(2\times n)$, Discrete Math.
281 (2004), 67--81.

\bibitem{Calkin}N.J. Calkin, Factors of sums of powers of binomial coefficients,
Acta Arith. 86 (1998), 17--26.


\bibitem{CD}M. Chamberland and K. Dilcher,
Divisibility properties of a class of binomial sums, J. Number Theory 120 (2006), 349--371.

\bibitem{CC}X. Chen and W. Chu, Moments on Catalan numbers, J. Math. Anal. Appl. 349 (2009), 311--316.

\bibitem{CH}W.Y.C. Chen and Q.-H. Hou, Factors of the Gaussian coefficients, Discrete Math. 306 (2006), 1446--1449.

\bibitem{FH}J. F\"urlinger and J. Hofbauer, $q$-Catalan numbers, J. Combin. Theory, Ser. A 2 (1985), 248--264.

\bibitem{GR}G. Gasper and M. Rahman, Basic Hypergeometric Series, Second Edition,
Encyclopedia of Mathematics and Its Applications, Vol. 96, Cambridge
University Press, Cambridge, 2004.

\bibitem{GJ}V.J.W. Guo and Q.-Q. Jiang, Factors of alternating sums of powers of $q$-Narayana numbers, J. Number Theory 177 (2017), 37--42.

\bibitem{GJZ}V.J.W. Guo, F. Jouhet, and J. Zeng, Factors of alternating sums
of products of binomial and $q$-binomial coefficients, Acta Arith. 127 (2007), 17--31.

\bibitem{GK}V.J.W. Guo and C. Krattenthaler, Some divisibility properties of binomial and $q$-binomial coefficients, J. Number Theory, 135 (2013) 167-184.

\bibitem{GZ04}V.J.W. Guo and J. Zeng, A $q$-analogue of Faulhaber's formula for sums of powers, Electron. J. Combin. 11 (2) (2004-2006), \#R19.

\bibitem{GJZ2005}V.J.W. Guo and J. Zeng, A combinatorial proof of a symmetric $q$-Pfaff-Saalschutz identity, Electron. J. Combin. 12 (2005), \#N2.

\bibitem{GZ2010}V.J.W. Guo and J. Zeng, Factors of binomial sums from the Catalan triangle, J. Number Theory 130 (2010), 172--186.

\bibitem{GZ06}V.J.W. Guo and J. Zeng, Some arithmetic properties of the $q$-Euler numbers and $q$-Sali\'e, numbers,
European J. Combin. 27 (2006), 884--895.

\bibitem{GHMR}J.M. Guti\'errez, M.A. Hern\'andez, P.J. Miana, N. Romero, New identities in the Catalan triangle,
J. Math. Anal. Appl. 341 (2008), 52--61.

\bibitem{KW}D. Knuth and H. Wilf, The power of a prime that divides a generalized binomial coefficient,
J. Reine Angew. Math. 396 (1989), 212--219.

\bibitem{MR}P.J. Miana and N. Romero, Computer proofs of new identities in the Catalan triangle,
Proceedings of the ``Segundas Jornadas de Teor\'ia de N\'umeros" (2007),
Bibl. Rev. Mat. Iberoamericana, Rev. Mat. Iberoamericana, Madrid, 2008, pp. 203--208.


\bibitem{Schlosser}M. Schlosser, $q$-Analogues of the sums of consecutive integers, squares, cubes, quarts and quints,
Electron. J. Combin. 11 (2004), \#R71.

\bibitem{Shapiro}L.W. Shapiro, A Catalan triangle, Discrete Math. 14 (1976), 83--90.

\bibitem{Sun0}Z.-W. Sun,  On Delannoy numbers and Schr\"oder numbers, J. Number Theory 131 (2011), 2387--2397.

\bibitem{Sun1}Z.-W. Sun, On sums of Ap\'ery polynomials and related congruences, J. Number Theory, 132 (2012),
2673--2699.

\bibitem{Sun2}Z.-W. Sun, Connections between $p=x^2+3y^2$ and Franel numbers, J. Number Theory 133 (2013), 2914--2928.

\bibitem{Sun3}Z.-W. Sun, Supercongruences involving products of two binomial coefficients, Finite Fields Appl. 22 (2013), 24-44.

\bibitem{Sun4}Z.-W. Sun, Congruences for Franel numbers, Adv. Appl. Math. 51 (2013),  524--535.

\bibitem{Sun5}Z.-W. Sun, Congruences involving generalized trinomial coefficients, Sci. China 57 (2014), 1375--1400.

\bibitem{WZ}S.O. Warnaar and W. Zudilin, A $q$-rious positivity, Aequationes Math. 81  (2011), 177--183.

\bibitem{Zeilberger}D. Zeilberger, A $q$-Foata proof of the $q$-Saalsch\"utz identity,
European J. Combin. 8 (1987), 461--463.

\bibitem{Zu}W. Zudilin, On a combinatorial problem of Asmus Schmidt, Electron. J. Combin. 11 (2004), \#R22.

\end{thebibliography}
\end{document}